\documentclass[11pt,a4paper,reqno]{amsart}

\usepackage{mathrsfs,amssymb}

\usepackage{pstricks,pst-node,pst-plot}
\usepackage{epic}
\usepackage{autograph}

\newtheorem{theorem}{Theorem}
\newtheorem{lemma}[theorem]{Lemma}
\newtheorem{proposition}[theorem]{Proposition}
\newtheorem{corollary}[theorem]{Corollary}

\theoremstyle{definition}

\newtheorem{definition}[theorem]{Definition}

\theoremstyle{plain}
\newcommand{\BZ}{\mathbb{Z}}
\newcommand{\ol}[1]{{\overline{#1}}}
\newcommand{\BN}{\mathbb{N}}

\newcommand{\BF}{\mathbb{F}}

\begin{document}
    \title{Formal languages and groups as memory}

\maketitle

\begin{center}
    Mark Kambites \\

    \medskip

    School of Mathematics, \  University of Manchester \\
    Manchester M60 1QD, \  England

    \medskip

    \texttt{Mark.Kambites@manchester.ac.uk} \\

\end{center}

\begin{abstract}
We present an exposition of the theory of \textit{$M$-automata} and
\textit{$G$-automata},
or finite automata augmented with a multiply-only register storing an element of a given
monoid or group. Included are a number of new results of a foundational
nature. We illustrate our techniques with a group-theoretic interpretation
and proof of a key theorem of Chomsky and Sch\"utzenberger from formal language
theory. 
\end{abstract}

\section{Introduction}

In recent years, both computer scientists and pure mathematicians have become
increasingly interested in the class of \textit{$M$-automata}, or finite state
automata augmented with a memory register which stores at any given time an
element of a given monoid
$M$. The register is initialised with the identity element of the
monoid; while reading an input word the automaton can modify the register
contents by multiplying by elements of the monoid. A word is accepted by the
automaton if, having read the entire word, the automaton reaches a final
state, with the register returned to the identity element.

Such automata have arisen repeatedly, both explicitly and implicitly,
in the theory of computation. For example, the blind 
$n$-counter machines studied by Greibach \cite{Greibach78} are simply
$\BZ^n$-automata. Related examples have also been studied by Ibarra,
Sahni and Kim \cite{Ibarra76}. $M$-automata are also equivalent to a
class of regulated grammars, known as
\textit{regular valence grammars} \cite{Fernau02}. More recently, there
has been increasing interest in this idea from
pure mathematicians, especially in the case that the register monoid is
a group. An area of lasting interest in combinatorial group theory is
the connection between structural properties of infinite discrete groups,
and language theoretic properties of their \textit{word problems}. Results of Gilman and Shapiro \cite{Gilman98}, of Elston and Ostheimer
\cite{Elston04} and of the author \cite{KambitesWordProblemsRecognisable}
have demonstrated that $M$-automata can play a useful role in this area.
Moreover, recent research of the author, Silva and Steinberg
\cite{KambitesGraphRat} has established a connection between the theory
of $G$-automata and the \textit{rational subset problem} for groups,
allowing language-theoretic results to be applied to decision problems
in group theory.

One aim of the present paper is to provide a self-contained introduction 
to the theory of this important area, in a form intelligible 
both to pure mathematicians and to computer scientists. In doing so, we 
aim to make explicit and transparent the connections between the computational
and the  algebraic approaches to $M$-automata and $G$-automata. This is
especially important since, to date, group theorists, semigroup theorists
and computer scientists have all worked on this topic, often unaware of
crucial results from the other discipline. Another objective is to establish
a number of new results of a foundational nature.

In addition to this introduction, this paper comprises three sections.
Section~\ref{sec_prelim} briefly recalls some necessary preliminaries,
and gives suggested references to more detailed treatments.
 Section~\ref{sec_basics} introduces $M$-automata and the families of languages they define. We also explain 
their relationship to the theory of rational transductions, and study the 
connection between algebraic properties of monoids and closure properties 
of the language classes they define. In Section~\ref{sec_contextfree}, we 
proceed to show how $M$-automata techniques can be used to obtain a 
group-theoretic interpretation and proof of one of the most important 
results in formal language theory -- namely, the Chomsky-Sch\"utzenberger 
theorem characterising context-free languages as the images under rational 
transductions of 2-sided Dyck languages \cite{Chomsky63}.
 
\section{Preliminaries}\label{sec_prelim}

In this section, we introduce the basic definitions which will be
required in this paper.  We begin with a very brief introduction to
formal languages and automata; a more comprehensive exposition can be found in
any of the numerous texts on the subject, such as \cite{Hopcroft69}. We
assume a familiarity with some basic definitions from algebra, such
as semigroups and monoids \cite{Howie95}, groups, generating sets and
presentations \cite{Lyndon77,Magnus76}.
More specialist notions from algebra will be defined as and when they are needed.

Let $\Sigma$ be a finite set of symbols, called an
\textit{alphabet}. A \textit{word} over $\Sigma$ is a finite sequence of
zero or more symbols from $\Sigma$; the unique \textit{empty word} of length
zero is denoted $\epsilon$. The set of all words over $\Sigma$
forms a monoid under the operation of \textit{concatenation}; this is
called the \textit{free monoid} on $\Sigma$ and denoted $\Sigma^*$.
A \textit{language} over $\Sigma$ is a set of words over $\Sigma$, that is,
a subset of the free monoid $\Sigma^*$.

A \textit{finite automaton} over a monoid $M$ is a finite directed
graph, possibly with loops and multiple edges, with each edge labelled
by element of $M$, together with a designated
\textit{initial vertex} and a set of designated \textit{terminal vertices}.
The vertices and edges of an automaton are often called \textit{states}
and \textit{transitions} respectively. The labelling of edges extends
naturally to a labelling of (directed) paths by elements of $M$. The \textit{subset accepted} or
\textit{recognised} by the automaton is the set of all elements of $M$
which label paths between the initial vertex and some terminal
vertex. A subset recognised by some automaton is called a
\textit{rational subset of $M$}. Notice that finitely generated submonoids
of $M$ are examples of rational subsets.

Of particular interest is the case where $M = \Sigma^*$ is a free monoid
on an alphabet $\Sigma$, so that the automaton accepts a set of words
over $\Sigma$, that is, a language over $\Sigma$. A language accepted by
such an automaton is called a \textit{rational language} or a
\textit{regular language}. The formal study of languages in
general, and of regular languages in particular, is of fundamental
importance in theoretical computer science, and increasingly also in
combinatorial algebra.

Another interesting case is that where $M = \Sigma^* \times \Omega^*$ is
a direct product of free monoids. Such an automaton is called
a \textit{finite transducer} from $\Sigma^*$ to $\Omega^*$; it recognises a relation, termed a
\textit{rational transduction}. Relations between free monoids, and
rational transductions in particular, are a
powerful tool for studying relationships between languages. If
$R \subseteq \Sigma^* \times \Omega^*$ then we say that the \textit{image}
of a language $L \subseteq \Sigma^*$ under $R$ is the language
of all words $v \in \Omega^*$ such that $(u, v) \in R$ for some $u \in L$.
We say that a language $K$ is a \textit{rational transduction} of a
language $L$ if $K$ is the image of $L$ under some rational transduction.
For a detailed exposition of
the theory of rational transductions, see \cite{Berstel79}.

\section{$M$-automata and $G$-automata}\label{sec_basics}

In this section, we introduce the definitions and some basic properties
of $M$-automata and $G$-automata. Let $M$ be a monoid with identity $1$
and $\Sigma$ a finite
alphabet. An \textit{$M$-automaton over $\Sigma$} is a finite
automaton over the direct product monoid
$M \times \Sigma^*$. For simplicity, we assume that the edges are labelled
by elements of $M \times (\Sigma \cup \lbrace \epsilon \rbrace)$.
We identify the free monoid $\Sigma^*$ with its natural embedding into
$M \times \Sigma^*$ as $\lbrace 1 \rbrace \times \Sigma^*$; thus, a word
$w \in \Sigma^*$ is \textit{accepted} by the automaton if there is a path
from the initial vertex to a terminal vertex labelled
$(1, w) \in M \times \Sigma^*$. The \textit{language accepted}
by the automaton is the set of all words in $\Sigma^*$ accepted by the
automaton; it is the intersection of the \textit{subset} accepted with
the embedded copy $\lbrace 1 \rbrace \times \Sigma^*$ of $\Sigma^*$.
We denote by $\BF(M)$ the family of all
languages accepted by $M$-automata.

From a mathematical perspective, then, the theory of $M$-automata can be 
viewed as an attempt to understand certain properties of a space (the free 
monoid $\Sigma^*$) by embedding it into a larger space ($M \times \Sigma^*$)
with more structure; in this sense, it is very roughly analoguous to the
embedding of the real numbers into the complex numbers. From a computational
perspective, an $M$-automaton can be thought of as a (non-deterministic)
finite automaton augmented with an extra memory register, which stores at any
point an element of the monoid $M$. The register is initialised with the
identity element of the monoid, and at each stage in its operation, the
automaton can modify the contents of the register by multiplication on
the right by some 
element of the monoid $M$. Of course it can also leave the register 
unchanged, simply by multiplying by the identity. The automaton cannot 
read the register during operation, but the contents act as an extra 
barrier to acceptance --- a word is accepted only if reading it can result 
in reaching a final state in which the register value has returned to the 
identity.

As an example, recall that 
a \textit{blind $n$-counter automaton} is a finite automaton
augmented with $n$ registers,
each of which stores a single integer value \cite{Greibach78}. The registers can be
incremented and decremented but not read; they are initialised to zero, and
a word is accepted exactly if, when reading it, the automaton can reach a
final state with all registers returned to zero. In view of the discussion
above, it is clear that a blind $n$-counter automaton is essentially the
same thing as a $\BZ^n$-automaton, where $\BZ$ denotes the group of integers
under addition, that is, the infinite cyclic group.

Notice that we have not required that the register monoid be finitely
generated. However, the following elementary observation will often
allow us to restrict attention to the case in which it is.
\begin{proposition}\label{prop_fingensub}
Let $M$ be a monoid, and suppose that $L$ is accepted by an $M$-automaton. Then there
exists a finitely generated submonoid $N$ of $M$ such that $L$ is accepted
by a $N$-automaton.
\end{proposition}
\begin{proof}
Since an $M$-automaton has finitely many edges, only finitely many elements
of $M$ can feature as the left-hand component of edge labels in a given
automaton. Clearly, the register can only ever hold values in the submonoid
$N$ of $M$ generated by these elements, so it suffices to view the automaton
as an $N$-automaton.
\end{proof}
Notwithstanding Proposition~\ref{prop_fingensub}, it is occasionally useful
to consider the \textit{class} of all $M$-automata where $M$ is not
finitely generated, since the corresponding class of languages (the union
of the classes corresponding to the finitely generated submonoids of $M$)
may not be defined by a single finitely generated monoid. See, for example,
\cite[Theorem~6.2]{Gilman98}, for an application of this approach.

If $M$ is a monoid generated by a set $X$, we say that the \textit{identity
language} $W_X(M)$ of $M$ with respect to $X$ is the set of all words over
$X$ representing the identity. In the case $M$ is a group, the identity
language is traditionally called
the \textit{word problem} of the group; this terminology is justified by
the fact that the membership problem of this language is algorithmically
equivalent to the problem of deciding whether two given words represent
the same element of the group, that is, to the word problem in the sense of
universal algebra. In a general monoid there is no such equivalence, and so
the term word problem is less appropriate.

The following simple observation has been made by several authors
(see, for example, \cite{Gilman98}), but apparently overlooked by a
number of others. It allows us to apply many standard results
from the theory of formal languages to the study of $M$-automata.
\begin{proposition}\label{prop_mauto_transduction}
Let $L$ be a language and $M$ a finitely generated monoid. Then the
following are equivalent:
\begin{itemize}
\item[(i)] $L$ is accepted by an $M$-automaton;
\item[(ii)] $L$ is a rational transduction of the identity language of $M$
            with respect to some finite generating set;
\item[(iii)] $L$ is a rational transduction of the identity language of $M$
            with respect to every finite generating set.
\end{itemize}
\end{proposition}
\begin{proof}
We begin by proving the equivalence of (i) and (ii).
First suppose (i) holds, and let $A$ be an $M$-automaton accepting $L$.
Let $S$ be the (necessarily finite) set of elements of $M$ which occur
on the left-hand-side of edge labels in $A$. Extend $S$ to a finite
generating set $X$ for $M$. Now $A$ can be viewed as a finite
transducer from $X^*$ to $\Sigma^*$. It follows easily from the
relevant definitions that the image of $W_X(M)$ is exactly
the language $L$, so that (ii) holds.

Conversely, suppose (ii) holds. Then there is a finite generating set $X$
for $M$ and a finite transducer $A$ from $X^*$ to $\Sigma^*$ such
that $L$ is the image of $W_X(M)$ under $A$. We obtain from $A$ an
$M$-automaton, by replacing each edge label $(w, x)$ with $(m,x)$ where
$m$ is the element of $M$ represented by $w \in X^*$. Again, it
follows easily from the definitions that the resulting $M$-automaton
accepts exactly the language $L$.

Since the monoid is assumed to be finitely generated, it is immediate
that (iii) implies (ii). It remains only to show that (ii) implies (iii).
Let $Y$ and $X$ be finite generating sets for a monoid $M$. For
each symbol $a \in X$, choose a word $w_a \in Y^*$ representing
the same element of $M$. Now let $R$ be the submonoid of
$Y^* \times X^*$
generated by the (finitely many) pairs of the form $(w_a, a)$. $R$ is
a finitely generated submonoid, and hence also, by our observations above, 
a rational transduction. Now if $L$ is the image of $W_X(M)$ under
a rational transduction $S$ then the relational composition
$$R \circ S = \lbrace (u, w) \mid (u, v) \in R, (v, w) \in S \text{ for some } v \in X^* \rbrace$$
is a rational transduction \cite[Theorem~II.4.4]{Berstel79}, and it is easily
verified that $L$ is the image of $W_Y(M)$ under $R \circ S$.
\end{proof}
Proposition~\ref{prop_mauto_transduction} tells us that the theory of
$M$-automata can be viewed as a special case  of the well-established field
of rational transductions. Indeed, we can easily and profitably translate
a large body of existing theory concerning rational transductions into the
$M$-automaton setting. For two main reasons, however, the study of
$M$-automata retains a distinct flavour, and remains of interest in its
own right. Firstly, the structure of the register monoid can be used to
prove interesting things about the accepting power of $M$-automata. Secondly,
$M$-automata can be used to gain insight into computational and language-theoretic
aspects of monoids. Both of these factors have special weight in the case
that the register monoid is a group, with all the extra structure that entails.

The following result, which has been observed independently by several authors
\cite{Corson05, Elder05}, is an immediate corollary of Proposition~\ref{prop_mauto_transduction} together with the fact that
rational transductions are closed under composition \cite[Theorem~III.4.4]{Berstel79}.
\begin{corollary}\label{cor_wordproball}
Let $M$ and $N$ be finitely generated monoids. Then the identity language of $N$ is accepted by
an $M$-automaton if and only if every language accepted by an $N$-automaton
is accepted by an $M$-automaton.
\end{corollary}

We end this section with a discussion of closure properties of language
families of the form $\BF(M)$. First, we observe that
every language class of the form $\BF(M)$ is easily seen to be closed under
finite union, as a simple consequence of non-determinism. On the other hand,
a class $\BF(M)$ need not be closed under intersection. The following
theorem will allow us to provide a straightforward characterisation of when
such a class is intersection-closed. The converse part was essentially
proved by Mitrana and Stiebe \cite{Mitrana01} in the case that the register
monoids are groups. Recall that a morphism between free monoids is called
\textit{alphabetic} if it maps each letter of the domain alphabet to
either a single letter or the empty word; one language is said to be
\textit{an alphabetic morphism of} another if the former is the image
of the latter under an alphabetic morphism of free monoids.

\begin{theorem}\label{thm_directproduct}
Let $M_1, M_2, \dots M_n$ be monoids. Then a language is accepted by an
$(M_1 \times \dots \times M_n)$-automaton if and only if it is an
alphabetic morphism of a language of the form $L_1 \cap \dots \cap L_n$
where each $L_i$ is accepted by an $M_i$-automaton.
\end{theorem}
\begin{proof}
Suppose first that $L \subseteq \Sigma^*$ is accepted by an
$(M_1 \times \dots \times M_n)$-automaton. It follows easily from
Proposition~\ref{prop_fingensub} that $L$ is
accepted by an $(M_1' \times \dots \times M_n')$-automaton where each
$M_i'$ is a finitely generated submonoid of $M_i$. For each $M_i'$
choose a finite generating set $X_i$ and let $X = X_1 \cup \dots \cup X_n$
so that $X$ is a finite generating set for $M_1' \times \dots \times M_n'$.
Now by Proposition~\ref{prop_mauto_transduction}, $L$ is the image of
the identity
language $W = W_X(M_1' \times \dots \times M_n')$ under some rational
transduction $\rho \subseteq X^* \times \Sigma^*$.

We claim that $W$ can be written as
$$W = K_1 \cap \dots \cap K_n$$
where each $K_i \in X^*$ lies in $\BF(M_i')$. Let $K_i$ be the set of
all words $w \in X^*$ such that $w$ represents a tuple in
$M_1' \times \dots \times M_n'$ with the identity element $1$ in the $i$th position.
It is readily verified that $K_i$ is accepted by an $M_i'$-automaton
with a single (initial and terminal) state $q$, and a loop at $q$
labelled $(x, g)$ whenever the letter $x \in X$ represents a tuple
with $g$ in the $i$th position. Now a word $w$ represents the identity of
$M_1' \times \dots \times M_n'$ if and only if $w$ lies in $K_i$ for all
$i$, so that $W$ is the intersection of the $K_i$ as required.

Now by \cite[Theorem~III.4.1]{Berstel79}, there exists an alphabet $Z$, a
regular language $R \subseteq Z^*$ and two alphabetic morphisms 
$\alpha : Z^* \to X^*$ and $\beta : Z^* \to \Sigma^*$ such
that
$$L \ = \ W \rho \ = \ (W \alpha^{-1} \cap R) \beta \ = \ ((K_1 \cap \dots \cap K_n) \alpha^{-1} \cap R) \beta.$$
It is easy to check that inverse morphisms distribute over intersection,
so that
$$(K_1 \cap \dots \cap K_n) \alpha^{-1} \ = \ K_1 \alpha^{-1} \cap \dots \cap K_n \alpha^{-1}.$$
Thus we obtain
$$L \ = \ ((K_1 \alpha^{-1} \cap \dots \cap K_n \alpha^{-1}) \cap R) \beta \ = \ ((K_1 \alpha^{-1} \cap R) \cap \dots \cap (K_n \alpha^{-1} \cap R)) \beta.$$
Each $\BF(M_i')$ is closed under rational transductions, and hence also under
inverse morphisms and intersection with regular languages.
It follows that
$$(K_i \alpha^{-1} \cap R) \ \in \ \BF(M_i') \ \subseteq \ \BF(M_i)$$
for each $i$, so setting $L_i = (K_i \alpha^{-1} \cap R)$
completes the proof of the direct implication.

For the converse, it is easy to check that the graph of an alphabetic
morphism between free monoids is a ratonal transduction; it follows that
the family $\BF(M_1 \times \dots \times M_n)$ is closed under alphabetic
morphisms. Thus, it suffices to suppose $L = L_1 \cap \dots \cap L_n$ where each $L_i$ is
in $\BF(M_i)$, and show that $L \in \BF(M_1 \times \dots \times M_n)$.
For this, it is clearly sufficient to prove the case $n=2$ and then
apply induction. Suppose, then, that $L_1, L_2 \subseteq \Sigma^*$ and
that each $L_i$ is accepted and by an $M_i$-automaton $A_i$ with state
set $Q_i$. Recall that the edges of each $A_i$ are labelled
by elements of $M_i' \times (\Sigma \cup \lbrace \epsilon \rbrace)$.
We define an $(M_1 \times M_2)$-automaton $A$ with:
\begin{itemize}
\item state set $Q_1 \times Q_2$;
\item an edge from $(s,u)$ to $(t,v)$ labelled $((x,y), a)$ whenever
      $A_1$ has an edge from $s$ to $t$ labelled $(x,a)$ and $A_2$ has
      an edge from $u$ to $v$ labelled $(y,a)$ for some letter $a \in \Sigma$;
\item an edge from $(s,u)$ to $(t,u)$ labelled $((x,1), \epsilon)$ whenever
      $A_1$ has an edge from $s$ to $t$ labelled $(x, \epsilon)$;
\item an edge from $(s,u)$ to $(s,v)$ labelled $((1,y), \epsilon)$ whenever
      $A_2$ has an edge from $u$ to $v$ labelled $(y, \epsilon)$;
\item start state $(p_0, q_0)$ where $p_0$ and $q_0$ are the start states
      of $A_1$ and $A_2$ respectively; and
\item final states $(p,q)$ such that $p$ and $q$ are final states of $A_1$
      and $A_2$ respectively.
\end{itemize}
It is an easy exercise to verify that $A$ accepts exactly the intersection
$L = L_1 \cap L_2$.
\end{proof}

Returning to our example, an elementary application of
Theorem~\ref{thm_directproduct} is the following
characterisation of the classes of languages accepted by blind
$n$-counter automaton.

\begin{corollary}\label{cor_greibach}
A language is recognised by a blind $n$-counter automaton if and only if
it is an alphabetic morphism of the intersection of $n$ languages recognised
by blind $1$-counter automata.
\end{corollary}

We also obtain a characterisation of those monoids $M$ for which
$\BF(M)$ is closed under intersection.

\begin{corollary}\label{cor_intersectionclosure}
Let $M$ be a finitely generated monoid. Then $\BF(M)$ is closed under finite intersection if
and only if there exists an $M$-automaton accepting the identity language
of $M \times M$.
\end{corollary}
\begin{proof}
If there exists an $M$-automaton accepting the identity language of
$M \times M$ then by Corollary~\ref{cor_wordproball}, every
language recognised by an $(M \times M)$-automaton is recognised by
an $M$-automaton. But by Theorem~\ref{thm_directproduct}, this
means that every intersection of two $M$-automaton languages is recognised
by an $M$-automaton. It follows by induction that $\BF(M)$ is closed under
finite intersection.

Conversely, the identity language of $M \times M$ is certainly recognised by
an $(M \times M)$-automaton, and hence by Theorem~\ref{thm_directproduct}
is an alphabetic morphism of the intersection of two languages in $\BF(M)$. Thus, if the latter is
intersection-closed then it contains the identity language of $M \times M$,
as required.
\end{proof}
Corollary~\ref{cor_intersectionclosure} has a particularly interesting
interpretation in the case that the monoid $M$ is a free group. We shall
see in Section~\ref{sec_contextfree} below that a language is context-free
if and only if it is recognised by a free group automaton. A well-known
theorem of Muller and Schupp \cite{Muller83}, combined with a subsequent
result of Dunwoody \cite{Dunwoody85}, tells us that a finitely generated
group has context-free word problem if and only if it is \textit{virtually free},
that is, has a free subgroup of finite index.
It follows that the fact that context-free languages are not intersection
closed can be viewed as a manifestation of the fact that a direct product
of virtually free groups is not, in general, virtually free.

\section{Free Groups and Context-free Languages}\label{sec_contextfree}

In this section we consider $M$-automata where $M$ is drawn from two 
particularly significant classes of monoids; namely, polycyclic monoids 
and free groups. From the perspective of algebra, these may be 
considered the motivating examples for the subject. We observe that an 
important theorem of Chomsky and Sch\"utzenberger \cite{Chomsky63} has a 
natural interpretation in terms of $M$-automata, and show how $M$-automata 
techniques can be used to provide an algebraic and automata-theoretic proof of 
the theorem.

We begin by recalling a basic definition
from automata theory. Let $X$ be a finite alphabet. A \textit{pushdown store}
or \textit{stack} with alphabet $X$ is a storage device which
stores, at any one time, a finite but unbounded sequence of symbols
from $X$, The basic operations permitted
are appending a new symbol to the right-hand end of the sequence
(``pushing'' a symbol), removing the rightmost symbol from the sequence
(''popping'') and reading the rightmost symbol on the stack.

The possible configurations of a pushdown store are naturally modelled by
elements of the free monoid $X^*$.
The operations of pushing and popping can be modelled by partial
functions, defined upon subsets
of $X^*$. Specifically, for each symbol $x \in X$ we define a
function
$$P_x : X^* \to X^* \ \ \ \ \ w \mapsto wx$$
which models the operation pushing the symbol $x$ onto the stack. The
corresponding operation of popping $x$ can only be performed when the
stack is in certain configurations -- namely, when it has an $x$ as the
rightmost symbol --
and so is modelled by a partial function.
$$Q_x : X^* x \to X^* \ \ \ \ \ wx \mapsto w.$$
The set of functions
$$\lbrace P_x, Q_x \mid x \in X \rbrace$$
generates a submonoid of the monoid of all partial functions on
$X^*$, under the natural operation of composition.
This monoid, which was first explicitly studied by Nivat and Perrot
\cite{Nivat70}, is called the \textit{polycyclic monoid on $X$} and
denoted $P(X)$. The elements of
$P(X)$ encapsulate the various sequences of operations which can
be performed upon a pushdown store with alphabet $X$.
The \textit{rank} of $P(X)$ is defined to be the size $|X|$
of the alphabet $X$; a polycyclic monoid is uniquely determined
(up to isomorphism) by its rank. Polycyclic monoids also arise naturally in the structural theory of
semigroups; of particular importance is that of rank $1$, which is
known as the \textit{bicyclic monoid}. For more general information see
\cite[Section 9.3]{Lawson98}.

There is a natural embedding of the free monoid $X^*$ into the
polycyclic monoid $P(X)$, which takes each symbol $x$ to $P_x$.
With this in mind, we shall identify $P_x$ with $x$ itself. The element
$Q_x$ is an \textit{inverse} to $P_x$ in the sense of inverse
semigroup theory \cite{Lawson98}; hence, we shall denote it
$x^{-1}$. Thus, the monoid $P(X)$ is simply generated by the set
$$\ol{X} = \lbrace x, x^{-1} \mid x \in X \rbrace.$$
and we shall have no further need of the notation $P_x$ and $Q_x$.

The inversion operation extends to the whole of $P(X)$, by defining
$(x^{-1})^{-1} = x$ for each $x \in X$, and $(x_1 \dots x_n)^{-1} = x_n^{-1} \dots x_1^{-1}$
for $x_1, \dots, x_n \in \ol{X}$. If $x$ and
$y$ are distinct elements of $X$ then the product $x y^{-1}$ is the
\textit{empty partial function} and forms a \textit{zero} in the
semigroup $P(X)$.

We shall define a \textit{pushdown automaton} to be a
polycyclic monoid automaton. The equivalence of this definition to the
standard one is straightforward and well-documented, for example in
\cite{Gilman98}. For those familiar with the usual definition, we remark
that the top symbol on the stack corresponds to the rightmost position of the
word, and that the automaton accepts with empty stack.

The attentive reader may have noticed that our polycyclic monoid model
of the pushdown store does not provide explicitly for ``reading'' the
contents of the stack. However, a polycyclic monoid automaton can use non-determinism
to test the rightmost stack symbol, by attempting
to pop every possible symbol and moving to different states
depending upon which succeeds; all but one attempt will result in the
register containing a zero value, which effectively constitutes failure.
This is an example of a more general phenomenon, in which the apparent
blindness of an $M$-automaton can be overcome by the use of
non-determinism.

The languages accepted by pushdown automata are called \textit{context-free}.
The class of context-free languages, which also admits an equivalent
definition in terms of generating grammars \cite{Berstel79}, 
is one of the most important languages classes in computer science.

We recall also a key notion from combinatorial group theory. Recall that
the \textit{free group} on the alphabet $X$ is the group defined by the
monoid presentation
$$\langle \ol{X} \mid xx^{-1} = x^{-1} x = 1 \text{ for all }  x \in X \rangle.$$
Generators from $X$ we shall call \textit{positive generators},
while those from $\ol{X} \setminus X$ are \textit{negative generators}.
Free groups are of central importance in combinatorial and geometric
group theory; see \cite{Lyndon77} or \cite{Magnus76} for a detailed
introduction.

The identity languages of the free group $F(X)$ and the polycyclic monoid
$P(X)$, with respect to the standard generating set $\ol{X}$, are
well-known in formal language theory. They are called
the \textit{2-sided Dyck language} (or just \textit{Dyck language}) on $\ol{X}$ and
the \textit{1-sided Dyck language} (or \textit{restricted Dyck language}
or \textit{semi-Dyck language}) on $\ol{X}$ respectively. The former
consists of all words over $\ol{X}$ which can be reduced to the empty
word by successive deletion of factors of the form $x x^{-1}$ or
$x^{-1} x$ where $x \in X$. The latter contains all words which can be
reduced to the empty word by deleting only factors of the form $x x^{-1}$
with $x \in X$. In particular, we see that the latter is (strictly)
contained in the former. Thus, any word over $\ol{X}$ representing the
identity in $P(X)$ also represents the identity in $F(X)$; the converse
does not hold.

A well-known theorem of Chomsky and Sch\"utzenberger states that the
context-free languages are exactly
the  rational transductions of 1-sided Dyck languages, and of 2-sided
Dyck languages \cite[Proposition~2]{Chomsky63}. By
Proposition~\ref{prop_mauto_transduction} this result has the following
interpretation in the $M$-automaton setting.
\begin{theorem}[Chomsky-Sch\"utzenberger 1963]\label{thm_chomschut}
Let $L$ be a language. Then the following are equivalent:
\begin{itemize}
\item[(i)] $L$ is context-free;
\item[(ii)] $L$ is accepted by a polycyclic monoid [of rank $2$] automaton;
\item[(iii)] $L$ is accepted by a free group [of rank $2$] automaton.
\end{itemize}
\end{theorem}

We have already remarked that the equivalence of (i) and (ii) is the usual
equivalence of context-free grammars and pushdown automata, and so is
well-known. It is also straightforward to show that a pushdown automaton can
simulate a free group automaton, so that (iii) implies (ii); this is left as
an exercise for the interested reader. The restriction to polycyclic monoids
[respectively, free groups] of rank 2 is a simple consequence of the
well-known fact that every polycyclic monoid [free group] of countable rank
embeds into the polycyclic monoid [free group] of rank 2. What remains,
which is the real burden of the proof, is to show that (i) and/or (ii)
implies (iii).

The original proof of Chomsky
and Sch\"utzenberger starts with a context-free grammar, and produces from it
an appropriate rational transduction; an example of this approach can be
found in \cite{Berstel79}.
A direct group-theoretic proof of this result was claimed by Dassow and
Mitrana \cite{Dassow00}; however, their construction was fundamentally
flawed \cite{Corson05}. A correct algebraic proof has recently
been provided by Corson \cite{Corson05}, who exhibited a
free group automaton accepting the identity language of a polycyclic monoid
automaton. The authors of both \cite{Dassow00} and \cite{Corson05} appear
to have overlooked the equivalence of the statement to the theorem of
Chomsky and Sch\"utzenberger.

Theorem~\ref{thm_chomschut} is quite surprising, in view
of our comments above regarding the method used by a polycyclic monoid
automaton to read the rightmost symbol of the stack. A polycyclic monoid
automaton apparently makes fundamental use of its ability to ``fail'', by
reaching a zero configuration of the register monoid. Since a free group
has no zero, a free group automaton seems to have no such capability,
and appears to be ``blind'' in a much more fundamental way. However, it
transpires that a carefully constructed interplay between the finite
state control and the group register can achieve the desired ``failing''
effect.

We remark also upon an interesting corollary to Theorem~\ref{thm_chomschut}.
It is well known 
that every recursively enumerable language is a homomorphic image (and 
hence a rational transduction) of the intersection of two context-free 
languages. Hence, combining Theorems~\ref{thm_directproduct} and 
\ref{thm_chomschut}, we immediately obtain the following result, which
was first observed in the group case by Mitrana and Stiebe \cite{Mitrana01}. 

\begin{theorem}
Let $M$ be a free group of rank $2$ or more, or a polycyclic monoid of
rank $2$ or more. Then $\BF(M \times M)$ is the class of all recursively
enumerable languages.
\end{theorem}

In the rest of this section, we present an alternative group- and 
automata-theoretic proof of Theorem~\ref{thm_chomschut}. In particular, we show explicitly 
how a free group automaton can simulate the operation of a pushdown 
automaton. In the process, we also obtain some technical results relating 
polycyclic monoids to free groups, which may be of independent interest. 
We begin by introducing a construction of a free group automaton from a 
polycyclic monoid automaton, that is, a pushdown automaton.

Suppose $L \subseteq \Sigma^*$ is the language accepted by a
$P(X)$-automaton $A$ with state set $Q$, that is, by a pushdown automaton
with state set $Q$ and stack alphabet $X$. Let $\#$ be a new symbol not
in $X$, and let
$X^\#$ denote the alphabet $X \cup \lbrace \# \rbrace$. We construct from
$A$ an new finite automaton $A'$ with edges labelled by elements of
$(\ol{X^\#}^* \times \Sigma^*)$. It has:
\begin{itemize}
\item state set $Q' = Q_- \cup Q_+$ where
$$Q_+ = \lbrace q_+ \mid q \in Q \rbrace \text{ and } Q_- = \lbrace q_- \mid q \in Q \rbrace$$
are disjoint sets in bijective correspondence with $Q$;
\item start state $q_+$ where $q$ is the start state of $A$;
\item final states of the form $q_-$ where $q$ is a final state of $A$;
\item an edge from $p_+$ to $q_+$ labelled $(x\#,w)$
      whenever $A$ has an edge from $p$ to $q$ labelled 
      $(x,w)$ with $x$ a positive generator.
\item an edge from $p_-$ to $q_+$ labelled $(x'\#, w)$  whenever $A$ has
      an edge from $p$ to $q$ labelled by 
      $(x',w)$ with $x'$ a negative generator;
\item an edge from $p_+$ to $q_+$ labelled $(\epsilon,w)$
      whenever $A$ has an edge from $p$ to $q$ labelled 
      $(\epsilon,w)$;
\item for each $q \in Q$, an edge from $q_+$ to $q_-$ labelled
      $(\epsilon,\epsilon)$; and
\item for each $q \in Q$, a loop at state $q_-$ labelled $(\#^{-1}, \epsilon)$.
\end{itemize}

\begin{figure}
\begin{picture}(80,40)
\thicklines
\setloopdiam{10}
\Large
\letstate PDA=(10,20)    \drawinitialstate(PDA){} \drawrepeatedstate(PDA){}
\letstate Plus=(40,20)       \drawinitialstate(Plus){$+$}
\letstate Minus=(70,20)      \drawrepeatedstate(Minus){$-$}
\drawloop[t](PDA){$(x, a)$}
\drawloop[b](PDA){$(x^{-1}, b)$}

\drawloop[t](Plus){$(x \#, a)$}
\drawcurvedtrans(Minus,Plus){$(x^{-1} \#, b)$}
\drawloop[b](Minus){$(\#^{-1}, \epsilon)$}
\drawcurvedtrans(Plus,Minus){$(\epsilon, \epsilon)$}
\end{picture}
\caption{A pushdown automaton (left) and a free group or pushdown automaton
(right), both accepting
         the 1-sided Dyck language on a two-letter alphabet.}
\end{figure}
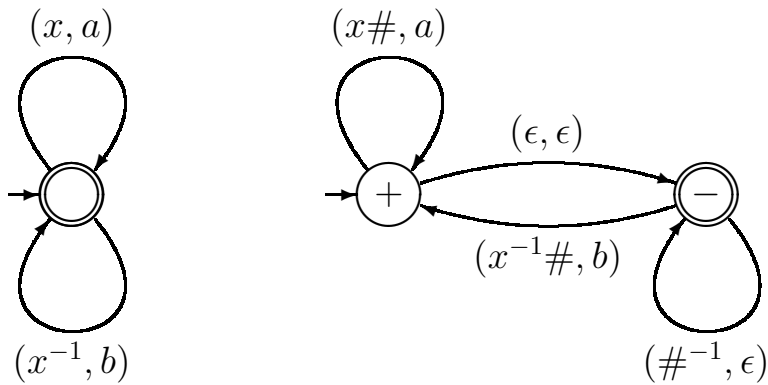
The automaton $A'$ can be interpreted either as a $P(X^\#)$-automaton or
as an $F(X^\#)$-automaton; it transpires that the language accepted
is the same for each choice.
Figure~1 illustrates a
$P(\lbrace x \rbrace)$-automaton accepting the 1-sided Dyck language on
the alphabet $\lbrace a, b \rbrace$, together with the automaton
constructed from it by the procedure above.

In general, we make the following claim.
\begin{theorem}\label{thm_constructionworks}
The free group automaton $A'$ and the polycyclic monoid automaton $A'$
both accept exactly the language $L$.
\end{theorem}
The rest of this section is devoted to the proof of
Theorem~\ref{thm_constructionworks}. We shall need a number of
preliminary definitions and results.
\begin{definition}
Let $w_1, \dots, w_n \in \ol{X}$ and $w = w_1 \dots w_n \in \ol{X}^*$.
Let $\#$ be a new symbol not in $\ol{X}$. A \textit{permissible padding}
of $w$ is a word of the form
$$x_1 x_2 \dots x_n (\#^{-1})^k$$
where $k \in \BN_0$ and for each $i \in \lbrace 1 \dots n \rbrace$ we have
$$x_i = \begin{cases}
          (\#^{-1})^m w_i \# \text{ for some } m \in \BN^0  &\text{ if $w_i$ is a negative generator;} \\
          w_i \#   &\text{ if $w_i$ is a positive generator.}
       \end{cases}
$$
\end{definition}
Thus, a permissible padding of $w$ is obtained by inserting the symbol
$\#$ after every generator in $w$, and zero or more $\#^{-1}$s before each
\textit{negative} generator, and at the end of the word.

The following lemma connects the above definition to our free group automaton
construction; it can be routinely verified.
\begin{lemma}\label{lemma_autopermissible}
The automaton $A'$ accepts $(x,w)$ if and only if there exists a
word $y \in \ol{X}^*$ such that $x$ is a permissible padding of $y$, and
$(y,w)$ is accepted by $A$.
\end{lemma}

We shall also use the following straightforward lemma concerning words
representing the identity in the free group.
\begin{lemma}\label{lemma_freegroupdecompose}
Let $w \in \ol{X}^*$ be a word representing the identity in a free group
$F(X)$, and
suppose $w = uxv$ where $u, v$ are words and $x \in \ol{X}$. Then either
$u$ has a suffix $x^{-1}e$ where $e$ represents the identity, or $v$ has
a prefix $ex^{-1}$ where $e$ represents the identity.
\end{lemma}
\begin{proof}
We have seen that any word in the 2-sided Dyck language, that is, any word
representing the identity in the free group, can be reduced to the empty
word by successively removing factors of the form $xx^{-1}$ and $x^{-1}x$
where $x \in X$. Such a reduction
process for $w$ must eventually bring the given occurrence of the generator
$x$ next to some occurrence $x^{-1}$, by deleting the letters between them.
But the product of these letters must be a factor representing the identity;
setting $e$ equal to this factor will give an appropriate factorisation of
either $u$ or $v$ (depending upon whether the given occurrence of $x^{-1}$
occurs before or after that of $x$).
\end{proof}

Recall that an element of the free group $F(X)$ is called \textit{positive}
if it can be written as a product of one or more positive generators. The
following definition facilitates a geometric interpretation of the positive
elements.
\begin{definition}
Let $w \in \ol{X}^*$ and let $x \in F(X)$. We say that $x$ is a \textit{minimum} of $w$, if
\begin{itemize}
\item[(i)] $w$ has a prefix representing $x$; and
\item[(ii)] no prefix of $w$ which represents $x$ is immediately followed
            by a negative generator.
\end{itemize}
\end{definition}

\begin{proposition}\label{prop_positive}
Let $w$ be a word representing the identity in $F(X)$. Then the
following are equivalent:
\begin{itemize}
\item[(i)] $w$ represents the identity in the polycyclic monoid $P(X)$;
\item[(ii)] every prefix of $w$ represents a positive or identity element;
\item[(iii)] the only minimum of $w$ is the identity of $F(X)$.
\end{itemize}
\end{proposition}
\begin{proof}
The equivalence of (i) and (ii) is well-known, and easily deduced from
the definitions.

Suppose now that (ii) holds, that is, that every prefix of $w$
represents a positive or identity element. It is easily seen that the
identity is a minimum of $w$. Moreover, if $x$ is a non-identity element
represented by a prefix of $w$ then consider the longest prefix of $w$
representing $x$. Considering the path traced through the Cayley graph
of $F(X)$, and recalling that the latter is a tree, it is clear that
the letter following this prefix must be a negative generator. Thus,
$x$ cannot be a minimum of $w$, and so (iii) is satisfied.

Conversely, suppose that (ii) does not hold, that is, that some prefix
of $w$ represents a non-positive, non-identity element. Suppose further
for a contradiction that (iii) holds, that is, that the identity is the
only minimum of $w$. Let $e$ be a non-positive, non-identity element
represented by a prefix of $w$. Since $e$ is not a minimum for $w$, there
is a prefix $u$ of $w$ representing $e$, which is followed by a negative
generator $x_1^{-1}$. But now $ux_1^{-1}$ represents another non-positive,
non-identity element. Continuing in this way, we obtain an infinite sequence
of prefixes of $u$, which must clearly all represent distinct elements. Since
$w$ is a finite word, this gives the required contradiction.
\end{proof}

\begin{lemma}\label{lemma_insertpad}
Let $w \in \ol{X^\#}^*$ be a word representing the identity in $P(X^\#)$
and suppose $w = uv$. Then there exists a factorisation $v = st$ such
that
\begin{itemize}
\item[(i)] either $t = \epsilon$ or $t$ begins with a negative generator; and
\item[(ii)]  $u \# s \#^{-1} t$ also represents the identity in $P(X^\#)$.
\end{itemize}
\end{lemma}
\begin{proof}
If $u$ represents the identity then $v$ also represents the identity,
so it suffices to take $s = v$ and $t = \epsilon$.

Assume now that $u$ does not represent the identity, and consider the path traced through the Cayley
graph of the free group $F(X^\#)$ when starting from the identity and
reading $w$. Since $w$ represents the identity in $P(X^\#)$, it also
represents the identity in $F(X^\#)$, so having reached the element
represented by $u$, this path must return to the identity. Since the
Cayley graph is a tree, the path must either leave in the
direction of the identity, in which case we take $s = \epsilon$, or leave
away from the identity and then return to the element represented by $u$
having read a word $s$, before leaving in the direction of the identity. By
Proposition~\ref{prop_positive}, $u$ represents a positive element,
so
``in the direction of the identity'' means following a negative
generator. Defining $t$ to be such that $v = st$, it is now clear
that $s$ and the corresponding $t$ have the desired properties.
\end{proof}
The following proposition, which may also be of interest in its own right,
is the main step in the proof.
\begin{proposition}\label{prop_identitypadding}
Let $w \in \ol{X}^*$. Then the following are equivalent.
\begin{itemize}
\item[(i)] $w$ represents $1$ in the polycyclic monoid $P(X)$;
\item[(ii)] $w$ admits a permissible padding which represents $1$ in the
polycyclic monoid $P(X^\#)$;
\item[(iii)] $w$ admits a permissible padding which represents $1$ in the
free group $F(X^\#)$.
\end{itemize}
\end{proposition}
\begin{proof}
First suppose (i) holds. By repeated application
of Lemma~\ref{lemma_insertpad}, we can insert the symbol $\#$ between every
pair of generators and the symbol $\#^{-1}$ in appropriate places, so as to
obtain a permissible padding of $w$ which represents $1$ in $P(X^\#)$. Thus,
(ii) holds.

Clearly every word representing the identity in $P(X^\#)$ also represents
the identity in $F(X^\#)$, so that (ii) implies (iii).

Finally, suppose (iii) holds, and let $w'$ be a permissible padding of $w$ which represents $1$ in
$F(X^\#)$. Suppose for a contradiction that $w$ does not represent
$1$ in $P(X)$. Certainly since $w'$ represents $1$ in $F(X^\#)$ we
must have that $w$ represents $1$ in $F(X)$. So by
Proposition~\ref{prop_positive}, $w$ contains a minimum which is
not the identity. Let $u$ be the shortest prefix of $w$ representing
this minimum, and write $w = uv$.

It follows that we can write $w' = u'\#v'$ where $u'$ and $v'$ are
paddings of $u$ and $v$ respectively. Certainly, since
$u$ is the shortest prefix representing the given element, $u$ has
no suffix representing the identity in $F(X)$. It follows that $u'$ has also
no suffix representing the identity in $F(X)$. Hence, by Lemma~\ref{lemma_freegroupdecompose}, we
can write $v' = e'\#^{-1}q'$ where $e'$ represents the identity.

Let $q$ and $e$ be the words over $\ol{X}^*$ obtained by deleting all
occurences of the letters $\#$ and $\#^{-1}$ from $q'$ and $e'$ respectively. Since $w' = u' \# e' \#^{-1} q'$
is a permissible padding of $w$, it follows that $q = \epsilon$ or $q$ begins
with a negative letter. But we have
$w = uv = ueq$ where $e$ represents
the identity. If $q = \epsilon$ then $u$ must represent the identity,
which is a contradiction. On the other hand, if $q$ begins with a negative
letter, then this contradicts the assumption that $u$ is a minimum
of $w$.
\end{proof}

We are now ready to complete the proof of Theorem~\ref{thm_constructionworks}.
\begin{proof}
Suppose a word $w$ is accepted by the pushdown automaton $A$. Then
by definition, there exists a word $x \in \ol{X}^*$ such that
$x$ represents the identity in $P(X)$, and $(x,w)$ is accepted by $A$
when viewed as a usual finite automaton. Now by Proposition~\ref{prop_identitypadding},
$x$ admits a permissible padding $y$ which represents $1$ in
the polycyclic monoid $P(X^\#)$, and hence also in the free group
$F(X^\#)$. Now by Lemma~\ref{lemma_autopermissible},
$(y,w)$ is accepted by $A'$ as a finite automaton over $F(X^\#) \times \Sigma^*$
and over $P(X^\#) \times \Sigma^*$. Hence, $w$ is accepted by $A'$ as both a
free group automaton and a polycyclic monoid automaton.

Conversely, if $w$ is accepted by $A'$ as a free group automaton 
[polycyclic monoid automaton], then by definition there exists a word $y$
such that $(y,w)$ is accepted by $A'$ as an automaton over
$F(X^\#) \times \Sigma^*$ [respectively, $P(X^\#) \times \Sigma^*$]
and $y$ represents $1$ in the free group [polycyclic monoid]. Now
by Lemma~\ref{lemma_autopermissible}, $y$ is a permissible padding of some
word $x$, such that $(x,w)$ is accepted by $A$ viewed as a finite
automaton over $P(X) \times \Sigma^*$. But by
Proposition~\ref{prop_identitypadding}, $x$ represents $1$ in the
polycyclic monoid $P(X)$, so that $w$ is accepted by $A$, as
required.
\end{proof}

\section*{Acknowledgements}

This research was conducted while the author was at Universit\"at Kassel,
and was supported by a Marie Curie Intra-European Fellowship within the
6th European Community Framework Programme. The author would like to thank
the many people with whom he had helpful discussions, and especially
Gretchen Ostheimer for her comments on the first draft of this article. He
would also like to thank Kirsty for all her support and encouragement.

\bibliographystyle{plain}

\def\cprime{$'$} \def\cprime{$'$}

\end{document}